# AN EXPLICIT FORMULA FOR THE SKOROKHOD MAP ON $[0, A]$

By Lukasz Kruk[1], John Lehoczky[2], Kavita Ramanan[3]
and Steven Shreve[4]

*Maria Curie-Sklodowska University and Carnegie Mellon University*

The Skorokhod map is a convenient tool for constructing solutions to stochastic differential equations with reflecting boundary conditions. In this work, an explicit formula for the Skorokhod map $\Gamma_{0,a}$ on $[0, a]$ for any $a > 0$ is derived. Specifically, it is shown that on the space $\mathcal{D}[0, \infty)$ of right-continuous functions with left limits taking values in $\mathbb{R}$, $\Gamma_{0,a} = \Lambda_a \circ \Gamma_0$, where $\Lambda_a : \mathcal{D}[0, \infty) \to \mathcal{D}[0, \infty)$ is defined by

$$\Lambda_a(\phi)(t) = \phi(t) - \sup_{s \in [0,t]} \left[ (\phi(s) - a)^+ \wedge \inf_{u \in [s,t]} \phi(u) \right]$$

and $\Gamma_0 : \mathcal{D}[0, \infty) \to \mathcal{D}[0, \infty)$ is the Skorokhod map on $[0, \infty)$, which is given explicitly by

$$\Gamma_0(\psi)(t) = \psi(t) + \sup_{s \in [0,t]} [-\psi(s)]^+.$$

In addition, properties of $\Lambda_a$ are developed and comparison properties of $\Gamma_{0,a}$ are established.

## 1. Introduction.

1.1. *Background.* In 1961 Skorokhod [13] considered the problem of constructing solutions to stochastic differential equations on the half-line $\mathbb{R}_+$

Received March 2006; revised June 2006.
[1]Supported by the State Committee for Scientific Research of Poland Grant 2 PO3A 012 23.
[2]Supported in part by ONR/DARPA MURI Contract N000140-01-1-0576.
[3]Supported in part by NSF Grants DMS-04-06191, DMI-03-23668-0000000965, DMS-04-05343.
[4]Supported in part by NSF Grant DMS-04-04682.

with a reflecting boundary condition at 0. His construction implicitly used properties of a deterministic mapping on the space $\mathcal{C}[0,\infty)$ of continuous functions on $[0,\infty)$. Anderson and Orey used this mapping more explicitly in their study of large deviations properties of reflected diffusions on a half-space in $\mathbb{R}^N$ (see page 194 of [1]). In particular, they exploited the fact that the mapping, which is now called the *Skorokhod map* and is denoted here by $\Gamma_0$, has the explicit representation

$$(1.1) \qquad \Gamma_0(\psi)(t) = \psi(t) + \sup_{s \in [0,t]} [-\psi(s)]^+, \qquad \psi \in \mathcal{C}[0,\infty),$$

and is consequently Lipschitz continuous (with constant 2) with respect to the uniform norm on $\mathcal{C}[0,\infty)$. El Karoui and Chaleyat-Maurel [6] used $\Gamma_0$ in a study of local times.

Given any trajectory $\psi$ in $\mathcal{D}[0,\infty)$, the space of right-continuous functions with left limits mapping $[0,\infty)$ into $\mathbb{R}$, $\Gamma_0$ can be extended using formula (1.1) to map $\psi$ to a "constrained version" $\phi = \psi + \eta$ of $\psi$ that is restricted to take values in $[0,\infty)$ by the minimal pushing term $\eta(t) \doteq \sup_{s \in [0,t]} [-\psi(s)]^+$. Minimality of $\eta$ follows from the fact that $\eta$ increases only at times $t$ when $\phi(t) = 0$ (see Definition 1.1 below for a precise statement). A multidimensional extension of the Skorokhod map was introduced by Tanaka [15]. Given any right-continuous function with left limits on $[0,\infty)$ taking values in $\mathbb{R}^N$, Tanaka produced a corresponding function taking values in a given convex domain by adding a constraining term on the boundary that acts in the direction normal to the boundary. Tanaka then used the solution to this Skorokhod problem to construct solutions of stochastic differential equations with normal reflection. In general, the Skorokhod map is a convenient tool for constructing processes that are restricted to take values in a certain domain by a constraining force that can push only along specified directions at the boundary. The study of many properties of the constrained or "reflected" process then reduces to the study of corresponding properties of the associated Skorokhod map.

In this paper, we focus on the particular case when the domain is a bounded interval in $\mathbb{R}$. For simplicity, for most of the paper we choose this interval to be $[0,a]$ for some $a > 0$, and denote the associated Skorokhod map by $\Gamma_{0,a}$. For functions in $\mathcal{D}[0,\infty)$, Chaleyat-Maurel, El Karoui and Marchal [4] posed and solved a version of this Skorokhod problem, producing functions taking values in $[0,a]$. However, in [4] the treatment of jumps across the boundary is different from that of Tanaka and this paper because in [4] the constrained function really "reflects" such jumps off the boundary, taking values in the interior of $[0,a]$, rather than being "constrained" to stay at the boundary. In contrast to [4], in this paper the Skorokhod map $\Gamma_{0,a}$ maps a trajectory in $\mathcal{D}[0,\infty)$ to a trajectory $\bar{\phi}$ in $\mathcal{D}[0,\infty)$ that is constrained to take values in $[0,a]$ by a *minimal* pushing force $\bar{\eta}$ that is allowed to increase only when $\bar{\phi}$ is at the lower boundary 0 and decrease only when $\bar{\phi}$ is at the upper boundary $a$ (see Definition 1.2 for a precise description of the Skorokhod map on $[0, a]$). Existence and uniqueness of solutions to this Skorokhod problem for continuous functions as well as step functions in $\mathcal{D}[0, \infty)$ follow directly from Lemmas 2.1, 2.3 and 2.6 of Tanaka [15]. In fact, it is well known that solutions to this Skorokhod problem exist for all functions in $\mathcal{D}[0, \infty)$ (see, e.g., [2]).

In this paper (see Theorem 1.4 and Remark 1.5 below) we provide an explicit formula for the Skorokhod map on a bounded interval in $\mathbb{R}$ (sometimes referred to as the *two-sided reflection map*). We first use this formula to provide direct proofs of Lipschitz continuity of this Skorokhod map and existence and uniqueness of solutions to the associated Skorokhod problem. In particular, our proofs do not rely on the existence and continuity results in [2] or [15], and also do not use approximation arguments. We then use this formula to establish comparison properties of $\Gamma_{0,a}$ (Theorem 1.7). This formula involves a new map, $\Lambda_a$, defined by (1.11). Properties of $\Lambda_a$ are developed in Proposition 1.3 and Corollary 1.6. In [5] a similar formula was obtained for the case when the "unconstrained" trajectory $\psi$ is of bounded variation. However, as elaborated in the next paragraph, in many applications of interest it is often important to understand the action of the two-sided reflection map on paths of unbounded variation.

The explicit formula for the Skorokhod map on $[0, \infty)$ has found application in a variety of contexts, including queueing theory and finance (see, e.g., [7, 8, 16]). More recently, it was used in [3] and [14] to derive various interesting distributional properties of quantities related to Brownian motion reflected on Brownian motion, a process that arises in the study of true self-repelling motions. In a similar fashion, the explicit formula for the Skorokhod map on a bounded interval in $\mathbb{R}$ is likely to have several potential applications. Already in [10] this formula plays a crucial role in the derivation of a diffusion approximation for the $GI/G/1$ queue with earliest-deadline-first service and reneging by customers who become late. In addition, in [9] the comparison properties of Theorem 1.7 are used to provide bounds on transaction costs in an optimal consumption/investment model. In the applications in both [9] and [10], the two-sided reflection map acts on paths of Brownian motion, which are almost surely of unbounded variation.

The outline of the paper is as follows. In Section 1.2 we introduce notation and recall the precise definitions and basic properties of $\Gamma_0$ and $\Gamma_{0,a}$. In Section 1.3, we state the main results. Properties of $\Lambda_a$ are established in Section 2. The proofs of Theorems 1.4 and 1.7 are presented in Sections 3 and 4 respectively. A technical result is relegated to the Appendix.

1.2. *Basic definitions.* Let $\mathcal{D}_+[0, \infty)$, $\mathcal{C}[0, \infty)$, $\mathcal{I}[0, \infty)$ and $\mathcal{BV}[0, \infty)$ denote the subspace of nonnegative, continuous, nondecreasing and bounded variation functions, respectively, in $\mathcal{D}[0,\infty)$. For $f \in \mathcal{BV}[0,\infty)$, $|f|_t$ denotes the total variation of $f$ on $[0,t]$. For $f \in \mathcal{D}[0,T]$, $\|f\|_T$ denotes the supremum norm of $f$ on $[0,T]$. Let $\mathbb{R}_+$ denote the set of nonnegative real numbers. Given $a, b \in \mathbb{R}$, denote $a \wedge b \doteq \min\{a,b\}$, $a \vee b \doteq \max\{a,b\}$, and $a^+ \doteq a \vee 0$. We denote by $\mathbb{I}_A$ the indicator function of a set $A$.

DEFINITION 1.1 (*Skorokhod map on* $[0,\infty)$). Given $\psi \in \mathcal{D}[0,\infty)$ there exists a unique pair of functions $(\phi, \eta) \in \mathcal{D}[0,\infty) \times \mathcal{I}[0,\infty)$ that satisfy the following two properties:

1. For every $t \in [0,\infty)$, $\phi(t) = \psi(t) + \eta(t) \in \mathbb{R}_+$;
2. $\eta(0-) = 0$, $\eta(0) \geq 0$, and

$$(1.2) \qquad \int_0^\infty \mathbb{I}_{\{\phi(s) > 0\}} \, d\eta(s) = 0.$$

The map $\Gamma_0 : \mathcal{D}[0,\infty) \to \mathcal{D}_+[0,\infty)$ that takes $\psi$ to the corresponding trajectory $\phi$ is referred to as the *one-sided reflection map* or *Skorokhod map on* $[0,\infty)$. The pair $(\phi, \eta)$ is said to solve the *Skorokhod problem on* $[0,\infty)$ for $\psi$.

Condition (1.2), often referred to as the *complementarity condition*, stipulates that the constraining term $\eta$ can increase only at times $t$ when $\phi(t) = 0$. As mentioned earlier, $\Gamma_0$, the Skorokhod map on $[0,\infty)$, has an explicit representation given by (1.1). The condition $\eta(0-) = 0$ is a convention by which we mean that $\eta(0) > 0$ implies that $\eta$ has a jump at zero and, according to (1.2), we must have $\phi(0) = 0$, in which case $\eta(0) = -\psi(0)$. This can happen only if $\psi(0) < 0$. In the event that $\psi(0) \geq 0$, we have $\eta(0) = 0$. In either case,

$$(1.3) \qquad \eta(0) = [-\psi(0)]^+.$$

In direct analogy with Definition 1.1 and the explicit representation (1.1) for $\Gamma_0$, it is easy to see that $\Gamma_a : \mathcal{D}[0,\infty) \to \mathcal{D}[0,\infty)$ defined by

$$(1.4) \qquad \Gamma_a(\psi)(t) \doteq \psi(t) - \sup_{s \in [0,t]} [\psi(s) - a]^+$$

takes $\psi \in \mathcal{D}[0,\infty)$ to the unique corresponding trajectory $\phi \in \mathcal{D}[0,\infty)$ that satisfies $\phi(t) \in (-\infty, a]$ for $t \in [0,\infty)$ and is such that $\eta = \psi - \phi$ is nondecreasing and increases only at times $t$ when $\phi(t) = a$ [i.e., such that $\int_0^\infty \mathbb{I}_{\{\phi(s) < a\}} \, d\eta(s) = 0$]. Indeed, it is straightforward to verify that given $a > 0$ and $\psi \in \mathcal{D}[0,\infty)$,

$$(1.5) \qquad \Gamma_a(\psi) = a - \Gamma_0(a - \psi).$$

The subject of this paper is the Skorokhod map that constrains a process in $\mathcal{D}[0,\infty)$ to remain within $[0,a]$, which is defined as follows.

DEFINITION 1.2 (*Skorokhod map $\Gamma_{0,a}$ on $[0,a]$*). Let $a > 0$ be given. Given $\psi \in \mathcal{D}[0,\infty)$ there exists a unique pair of functions $(\bar{\phi}, \bar{\eta}) \in \mathcal{D}[0,\infty) \times \mathcal{BV}[0,\infty)$ that satisfy the following two properties:

1. For every $t \in [0,\infty)$, $\bar{\phi}(t) = \psi(t) + \bar{\eta}(t) \in [0,a]$;
2. $\bar{\eta}(0-) = 0$ and $\bar{\eta}$ has the decomposition $\bar{\eta} = \bar{\eta}_\ell - \bar{\eta}_u$ as the difference of functions $\bar{\eta}_\ell, \bar{\eta}_u \in \mathcal{I}[0,\infty)$ satisfying

$$(1.6) \qquad \int_0^\infty \mathbb{I}_{\{\bar{\phi}(s)>0\}} \, d\bar{\eta}_\ell(s) = 0 \quad \text{and} \quad \int_0^\infty \mathbb{I}_{\{\bar{\phi}(s)<a\}} \, d\bar{\eta}_u(s) = 0.$$

We refer to the mapping $\Gamma_{0,a} : \mathcal{D}[0,\infty) \to \mathcal{D}[0,\infty)$ that takes $\psi$ to the corresponding $\bar{\phi}$ as the *two-sided reflection map* or the *Skorokhod map on $[0,a]$*. The pair $(\bar{\phi}, \bar{\eta})$ is said to solve the *Skorokhod problem on $[0,a]$* for $\psi$.

Similarly to (1.3), the condition $\bar{\eta}(0-) = 0$ coupled with the complementarity conditions (1.6) implies that

$$(1.7) \qquad \bar{\eta}(0) = [-\psi(0)]^+ - [\psi(0) - a]^+.$$

In other words, $\bar{\phi}(0) = \pi(\psi(0))$, where $\pi : \mathbb{R} \to [0,a]$ is the projection map

$$(1.8) \qquad \pi(x) = \begin{cases} a, & \text{if } x \geq a, \\ x, & \text{if } 0 \leq x \leq a, \\ 0, & \text{if } x \leq 0. \end{cases}$$

Furthermore, from the explicit expressions for $\Gamma_0$ and $\Gamma_a$ given in (1.1) and (1.4), respectively, it is clear (see, e.g., Section 2.3 of [7]) that $\bar{\eta}_\ell$ and $\bar{\eta}_u$ satisfy the equations

$$(1.9) \quad \bar{\eta}_\ell(t) = \sup_{s \in [0,t]} [\bar{\eta}_u(s) - \psi(s)]^+ \quad \text{and} \quad \bar{\eta}_u(t) = \sup_{s \in [0,t]} [\psi(s) + \bar{\eta}_\ell(s) - a]^+.$$

Now consider $\psi \in \mathcal{D}[0,\infty)$ and let $\bar{\eta} \doteq \Gamma_{0,a}(\psi) - \psi$, which has the decomposition $\bar{\eta} = \bar{\eta}_\ell - \bar{\eta}_u$ into the difference of processes in $\mathcal{I}[0,\infty)$ as in Definition 1.2. Denote $\tilde{\eta} \doteq \Gamma_{0,a}(a - \psi) - a + \psi$, which has the corresponding decomposition $\tilde{\eta} = \tilde{\eta}_\ell - \tilde{\eta}_u$. In a similar fashion to (1.5), it follows immediately from the definition that $\Gamma_{0,a}(\psi) = a - \Gamma_{0,a}(a - \psi)$ and, moreover, that

$$(1.10) \qquad \tilde{\eta}_\ell = \bar{\eta}_u \quad \text{and} \quad \tilde{\eta}_u = \bar{\eta}_\ell.$$

1.3. *Main results.* Our main result provides an explicit representation for the Skorokhod map $\Gamma_{0,a}$ on $[0,a]$ in terms of the mapping $\Lambda_a : \mathcal{D}[0,\infty) \to \mathcal{D}[0,\infty)$ defined by

$$(1.11) \qquad \Lambda_a(\phi)(t) \doteq \phi(t) - \sup_{s \in [0,t]} \left[ (\phi(s) - a)^+ \wedge \inf_{u \in [s,t]} \phi(u) \right].$$

For $t \in [0,\infty)$ and $s \in [0,t]$, we will use the notation

$$(1.12) \qquad R_t(\phi)(s) \doteq (\phi(s) - a)^+ \wedge \inf_{u \in [s,t]} \phi(u),$$

in terms of which (1.11) may be written as $\Lambda_a(\phi)(t) \doteq \phi(t) - \sup_{s \in [0,t]} R_t(\phi)(s)$. We list properties of $\Lambda_a$ and then state our main result as Theorem 1.4.

PROPOSITION 1.3. $\Lambda_a$ *maps* $\mathcal{D}[0,\infty)$ *into* $\mathcal{D}[0,\infty)$, $\mathcal{C}[0,\infty)$ *into* $\mathcal{C}[0,\infty)$, $\mathcal{BV}[0,\infty)$ *into* $\mathcal{BV}[0,\infty)$, *and absolutely continuous functions to absolutely continuous functions.*

The proof of Proposition 1.3 is the subject of Section 2.

THEOREM 1.4. *Given* $a > 0$, *let* $\Gamma_0$ *and* $\Gamma_{0,a}$ *be the Skorokhod maps on* $[0,\infty)$ *and* $[0,a]$ *respectively. Then*

$$(1.13) \qquad \Gamma_{0,a} = \Lambda_a \circ \Gamma_0.$$

REMARK 1.5. Consideration of the formula in [5] leads to a formula for $\Gamma_{0,a}$ different from (1.13) that can be derived from (1.13), namely (see [11])

$$\Gamma_{0,a}(\psi)(t) = \psi(t) - \left[(\psi(0) - a)^+ \wedge \inf_{u \in [0,t]} \psi(u)\right]$$
$$\vee \sup_{s \in [0,t]} \left[(\psi(s) - a) \wedge \inf_{u \in [s,t]} \psi(u)\right].$$

It is straightforward to generalize our results to the case where $[0,a]$ is replaced by $[z,a]$ for $-\infty < z < a < \infty$. In this case, the corresponding one-sided Skorokhod map $\Gamma_z$ is defined as in Definition 1.1, but with $\mathbb{R}_+$ replaced by $[z,\infty)$ in property 1 and $\phi(s) > 0$ replaced by $\phi(s) > z$ in equation (1.2), and the corresponding two-sided Skorokhod map $\Gamma_{z,a}$ is defined as in Definition 1.2, but with $[0,a]$ replaced by $[z,a]$ in property 1 and $\overline{\phi}(s) > 0$ replaced by $\overline{\phi}(s) > z$ in equation (1.6). A straightforward extension of Theorem 1.4 then shows that $\Gamma_{z,a}(\psi) = \Lambda_{z,a} \circ \Gamma_z(\psi)$, where $\Gamma_z$ and $\Lambda_{z,a}$ mapping $\mathcal{D}[0,\infty)$ into itself are defined by $\Gamma_z(\psi)(t) \doteq \psi(t) + \sup_{s \in [0,t]}[z - \psi(s)]^+$ and

$$(1.14) \qquad \Lambda_{z,a}(\phi)(t) \doteq \phi(t) - \sup_{s \in [0,t]} \left[(\phi(s) - a)^+ \wedge \inf_{u \in [s,t]} (\phi(u) - z)\right].$$

Theorem 1.4 allows us to give concise proofs of the Lipschitz continuity of the map $\Gamma_{0,a}$ in the uniform, $J_1$ and $M_1$ metrics. For this, we let $d_\infty$ denote the uniform metric on $[0,T]$, $d_0$ the standard $J_1$ metric on $\mathcal{D}[0,T]$ (see, e.g., definition (3.2) on page 79 of [16]), and $d_1$ the standard $M_1$ metric on $\mathcal{D}[0,T]$ (see, e.g., definition (3.4) on page 82 of [16]), while $\bar{d}_\infty$, $\bar{d}_0$ and $\bar{d}_1$ denote the corresponding metrics on $\mathcal{D}[0,\infty)$ (see, e.g., Section 12.9 of [16]).

COROLLARY 1.6. *There exists a constant $L$ such that for all $T > 0$ and $\psi_1, \psi_2 \in \mathcal{D}[0, T]$,*

(1.15) $\quad d_i(\Lambda_a(\psi_1), \Lambda_a(\psi_2)) \leq 2 d_i(\psi_1, \psi_2) \qquad$ *for $i = 0, 1, \infty$*;

(1.16) $\quad d_i(\Gamma_{0,a}(\psi_1), \Gamma_{0,a}(\psi_2)) \leq L d_i(\psi_1, \psi_2) \qquad$ *for $i = 0, 1, \infty$.*

*Moreover, the six inequalities above continue to hold for $\psi_1, \psi_2 \in \mathcal{D}[0, \infty)$ if $d_\infty$, $d_0$ and $d_1$ are replaced by $\bar{d}_\infty$, $\bar{d}_0$ and $\bar{d}_1$, respectively.*

The proofs of Theorem 1.4 and Corollary 1.6 are given in Section 3. Continuity of $\Gamma_0$ and $\Gamma_{0,a}$ in the $J_1$ and $M_1$ metrics is due to [2]. For proofs of the inequalities for $\Gamma_{0,a}$ in Corollary 1.6 that are different from the proofs in this paper, see Section 14.8 of [16].

Lastly, in Theorem 1.7, we state comparison properties of the Skorokhod map on $[0, a]$. The proof of this result is presented in Section 4.

THEOREM 1.7. *Given $a > 0$, $c_0, c_0' \in \mathbb{R}$ and $\psi, \psi' \in \mathcal{D}[0, \infty)$ with $\psi(0) = \psi'(0) = 0$, suppose $(\bar{\phi}, \bar{\eta})$ and $(\bar{\phi}', \bar{\eta}')$ solve the Skorokhod problem on $[0, a]$ for $c_0 + \psi$ and $c_0' + \psi'$, respectively. Moreover, suppose $\bar{\eta} = \bar{\eta}_\ell - \bar{\eta}_u$ is the decomposition of $\bar{\eta}$ into the difference of processes in $\mathcal{I}[0, \infty)$ satisfying (1.6) and $\bar{\eta}_\ell' - \bar{\eta}_u'$ is the corresponding decomposition of $\bar{\eta}'$. If there exists $\nu \in \mathcal{I}[0, \infty)$ such that $\psi = \psi' + \nu$, then the following four inequalities hold:*

1. $\bar{\eta}_\ell - [c_0' - c_0]^+ \leq \bar{\eta}_\ell' \leq \bar{\eta}_\ell + \nu + [c_0 - c_0']^+$;
2. $\bar{\eta}_u' - [c_0' - c_0]^+ \leq \bar{\eta}_u \leq \bar{\eta}_u' + \nu + [c_0 - c_0']^+$;
3. $\bar{\eta} - [c_0' - c_0]^+ \leq \bar{\eta}' \leq \bar{\eta} + \nu + [c_0 - c_0']^+$;
4. $[-[c_0 - c_0']^+ - \nu] \vee [-a] \leq \bar{\phi}' - \bar{\phi} \leq [c_0 - c_0']^+ \wedge a$.

**2. Proof of Proposition 1.3.** Let $\phi \in \mathcal{D}[0, \infty)$ be given. For each $\theta_1 \geq 0$ and $\varepsilon > 0$, there exists $\theta_2 > \theta_1$ such that

(2.1) $$\sup_{s, u \in [\theta_1, \theta_2)} |\phi(s) - \phi(u)| \leq \varepsilon.$$

Similarly, for each $\theta_2 > 0$ and $\varepsilon > 0$, there exists $\theta_1 \in [0, \theta_2)$ such that (2.1) holds. It is straightforward to use this observation and the following lemma to verify that $\Lambda_a(\phi)$ is right-continuous with left-hand limits, that is, that $\Lambda_a$ maps $\mathcal{D}[0, \infty)$ into $\mathcal{D}[0, \infty)$.

LEMMA 2.1. *Let $\phi \in \mathcal{D}[0, \infty)$ be given. For any $0 \leq \theta_1 < \theta_2$,*

$$\sup_{t_1, t_2 \in [\theta_1, \theta_2)} |\Lambda_a(\phi)(t_1) - \Lambda_a(\phi)(t_2)| \leq 2 \sup_{s, u \in [\theta_1, \theta_2)} |\phi(s) - \phi(u)|.$$

PROOF. From the definition of $R_t$ in (1.12), we see that for any $t \geq 0$,

$$(2.2) \qquad (\phi(t) - a)^+ \wedge \phi(t) \leq \sup_{s \in [0,t]} R_t(\phi)(s) \leq \phi(t).$$

Let $\varepsilon \doteq \sup_{s,u \in [\theta_1, \theta_2)} |\phi(s) - \phi(u)|$ and let $t_1, t_2$ be in $[\theta_1, \theta_2)$ with $t_1 \leq t_2$. Then $R_{t_2}(\phi)(s) \leq R_{t_1}(\phi)(s)$ for $s \in [0, t_1]$ and

$$\sup_{s \in (t_1, t_2]} R_{t_2}(\phi)(s) \leq \sup_{s \in (t_1, t_2]} (\phi(s) - a)^+ \leq (\phi(t_1) - a)^+ + \varepsilon.$$

Therefore

$$\sup_{s \in [0,t_2]} R_{t_2}(\phi)(s) = \sup_{s \in [0,t_1]} R_{t_2}(\phi)(s) \vee \sup_{s \in (t_1, t_2]} R_{t_2}(\phi)(s)$$

$$\leq \sup_{s \in [0,t_1]} R_{t_1}(\phi)(s) \vee [(\phi(t_1) - a)^+ + \varepsilon]$$

$$\leq \sup_{s \in [0,t_1]} R_{t_1}(\phi)(s) + \varepsilon,$$

where the last inequality uses the first inequality in (2.2). In turn, this yields

$$\Lambda_a(\phi)(t_2) = \phi(t_2) - \sup_{s \in [0,t_2]} R_{t_2}(\phi)(s)$$

$$(2.3) \qquad \geq \phi(t_1) - \varepsilon - \sup_{s \in [0,t_1]} R_{t_1}(\phi)(s) - \varepsilon$$

$$= \Lambda_a(\phi)(t_1) - 2\varepsilon.$$

The second inequality in (2.2) and the definition of $\varepsilon$ imply that

$$\sup_{s \in [0,t_1]} R_{t_1}(\phi)(s) - \varepsilon \leq \sup_{s \in [0,t_1]} R_{t_1}(\phi)(s) \wedge (\phi(t_1) - \varepsilon)$$

$$\leq \sup_{s \in [0,t_1]} \left[ R_{t_1}(\phi)(s) \wedge \inf_{s \in (t_1, t_2]} \phi(u) \right]$$

$$= \sup_{s \in [0,t_1]} R_{t_2}(\phi)(s)$$

$$\leq \sup_{s \in [0,t_2]} R_{t_2}(\phi)(s).$$

From this we conclude that

$$\Lambda_a(\phi)(t_2) = \phi(t_2) - \sup_{s \in [0,t_2]} R_{t_2}(\phi)(s)$$

$$\leq \phi(t_1) + \varepsilon - \sup_{s \in [0,t_1]} R_{t_1}(\phi)(s) + \varepsilon$$

$$= \Lambda_a(\phi)(t_1) + 2\varepsilon.$$

Together with (2.3) and the definition of $\varepsilon$, this proves the lemma. $\square$

REMARK 2.2. The proof of Lemma 2.1 also shows that the oscillation of $\Lambda_a(\phi)$ is bounded by the oscillation of $\phi$ on the closed interval $[\theta_1, \theta_2]$, that is,

$$\sup_{t_1,t_2\in[\theta_1,\theta_2]} |\Lambda_a(\phi)(t_1) - \Lambda_a(\phi)(t_2)| \leq 2 \sup_{s,u\in[\theta_1,\theta_2]} |\phi(s) - \phi(u)|.$$

Therefore $\Lambda_a$ maps $\mathcal{C}[0,\infty)$ to $\mathcal{C}[0,\infty)$.

COROLLARY 2.3. *$\Lambda_a$ maps absolutely continuous functions to absolutely continuous functions.*

PROOF. Suppose $\phi \in \mathcal{D}[0,\infty)$ is absolutely continuous. We fix an arbitrary $T > 0$. By the definition of absolute continuity, there exists a function $v_\phi : (0,\infty) \to (0,\infty)$ such that for every $\varepsilon > 0$ and every set of nonoverlapping intervals $(s_j, t_j)$, $j = 1, \ldots, J$, contained in $[0, T]$,

$$(2.4) \qquad \sum_{j=1}^J (t_j - s_j) < v_\phi(\varepsilon) \quad \implies \quad \sum_{j=1}^J |\phi(t_j) - \phi(s_j)| < \varepsilon.$$

Define the function $v_{\Lambda_a(\phi)} : (0,\infty) \to (0,\infty)$ by $v_{\Lambda_a(\phi)}(\varepsilon) \doteq v_\phi(\varepsilon/2)$ for $\varepsilon > 0$. We claim that (2.4) holds with $\phi$ replaced everywhere by $\Lambda_a(\phi)$, thus showing that $\Lambda_a(\phi)$ is absolutely continuous. For the proof of the claim, fix $\varepsilon > 0$ and consider any set of nonoverlapping intervals $(s_j, t_j)$, $j = 1, \ldots, J$, such that $\sum_{j=1}^J (t_j - s_j) < v_{\Lambda_a(\phi)}(\varepsilon)$. For $j = 1, \ldots, J$, choose $s_j \leq \overline{s}_j \leq \overline{t}_j \leq t_j$ such that $|\phi(\overline{t}_j) - \phi(\overline{s}_j)| = \max_{u,r \in [s_j, t_j]} |\phi(r) - \phi(u)|$. Remark 2.2 implies that $\Lambda_a(\phi) \in \mathcal{C}[0,\infty)$ and

$$\sum_{j=1}^J |\Lambda_a(\phi)(t_j) - \Lambda_a(\phi)(s_j)| \leq \sum_{j=1}^J \max_{u,r\in[s_j,t_j]} |\Lambda_a(\phi)(r) - \Lambda_a(\phi)(u)|$$

$$\leq 2 \sum_{j=1}^J \max_{u,r\in[s_j,t_j]} |\phi(r) - \phi(u)|$$

$$= 2 \sum_{j=1}^J |\phi(\overline{t}_j) - \phi(\overline{s}_j)|$$

$$\leq \varepsilon,$$

where the last inequality is a consequence of (2.4) and the fact that $\sum_{j=1}^J (\overline{t}_j - \overline{s}_j) < v_{\Lambda(\phi)}(\varepsilon) = v_\phi(\varepsilon/2)$. □

To complete the proof of Proposition 1.3, it remains only to show that $\Lambda_a$ maps $\mathcal{BV}[0,\infty)$ to $\mathcal{BV}[0,\infty)$. We do not use this fact in the present paper, and hence can use any results in the remainder of the paper to establish it. Recall the definition of $R_t(\phi)$ given in (1.12). For $\phi \in \mathcal{D}[0, \infty)$, it will be convenient to introduce the function $C^\phi \in \mathcal{D}[0, \infty)$ defined for $t \in [0, \infty)$ by

$$(2.5) \quad C^\phi(t) \doteq \sup_{s \in [0,t]} [R_t(\phi)(s)] = \sup_{s \in [0,t]} \left[ (\phi(s) - a)^+ \wedge \inf_{u \in [s,t]} \phi(u) \right].$$

Note that then $\Lambda_a(\phi) = \phi - C^\phi$ for every $\phi \in \mathcal{D}[0, \infty)$. According to Theorem 3.4 below, the function $C^\phi$ given by (2.5) has bounded variation. If $\phi$ also has bounded variation, then $\Lambda_a(\phi) = \phi - C^\phi$ does as well.

**3. Proof of Theorem 1.4.** An intuitive way of constructing $\bar{\phi} \doteq \Lambda_a(\phi)$ from $\phi \doteq \Gamma_0(\psi)$ is to first create two increasing sequences of times $\{\sigma_k\}_{k=0}^\infty$ and $\{\tau_k\}_{k=1}^\infty$ so that on each interval of the form $[\sigma_{k-1}, \tau_k)$, there is only pushing of $\phi$ from above and on each interval of the form $[\tau_k, \sigma_k)$, there is only pushing of $\phi$ from below. In this section we execute that construction and thereby obtain the decomposition in (3.24) below of the bounded variation process $C^\phi$ defined by (2.5) into the difference of two nondecreasing processes. For this construction, we assume that $\phi$ is in $\mathcal{D}_+[0, \infty)$. We have in mind that $\phi = \Gamma_0(\psi)$ for some $\psi \in \mathcal{D}[0, \infty)$.

For $\phi \in \mathcal{D}_+[0, \infty)$ and $a > 0$, we set $\tau_0 \doteq 0$,

$$(3.1) \quad \sigma_0 \doteq \min\{t \geq 0 | \phi(t) - a \geq 0\},$$

and for $k \geq 1$, we set

$$(3.2) \quad \tau_k \doteq \min\left\{ t \geq \sigma_{k-1} \Big| \phi(t) \leq \sup_{s \in [\sigma_{k-1}, t]} \phi(s) - a \right\},$$

$$(3.3) \quad \sigma_k \doteq \min\left\{ t \geq \tau_k \Big| \phi(t) - a \geq \inf_{u \in [\tau_k, t]} \phi(u) \right\}.$$

The minima in (3.1)–(3.3) over $t$ are obtained (or are $+\infty$) because of the right-continuity of $\phi$. In particular, for $k \geq 1$,

$$(3.4) \quad \sup_{s \in [\sigma_{k-1}, u]} \phi(s) - a < \phi(u) \quad \forall u \in [\sigma_{k-1}, \tau_k),$$

$$(3.5) \quad \sup_{s \in [\sigma_{k-1}, \tau_k]} \phi(s) - a \geq \phi(\tau_k),$$

$$(3.6) \quad \phi(s) - a < \inf_{u \in [\tau_k, s]} \phi(u) \quad \forall s \in [\tau_k, \sigma_k),$$

$$(3.7) \quad \phi(\sigma_k) - a \geq \inf_{u \in [\tau_k, \sigma_k]} \phi(u).$$

Furthermore,

$$(3.8) \quad \phi(\sigma_0) - a \geq 0.$$

We have $0 = \tau_0 \leq \sigma_0 < \tau_1 < \sigma_1 < \tau_2 < \sigma_2 < \cdots$.

PROPOSITION 3.1. *As $k \to \infty$, we have $\tau_k \uparrow \infty$ and $\sigma_k \uparrow \infty$.*

PROOF. Assume the proposition is false. Then there is a number $\theta < \infty$ such that $\tau_k \uparrow \theta$ and $\sigma_k \uparrow \theta$. Relation (3.5) implies the existence of $\rho_k \in [\sigma_{k-1}, \tau_k]$ such that $\phi(\rho_k) \geq \phi(\tau_k) + \frac{a}{2}$. Since $\rho_k \uparrow \theta$, $\phi$ does not have a left-hand limit at $\theta$. This contradicts the membership of $\phi$ in $\mathcal{D}_+[0, \infty)$. □

PROPOSITION 3.2. *For $k \geq 1$, $C^\phi(t) = \sup_{s \in [\sigma_{k-1}, t]}(\phi(s) - a)^+$ for all $t \in [\sigma_{k-1}, \tau_k)$.*

PROOF. Let $t \in (\sigma_{k-1}, \tau_k)$ and $\rho \in (\sigma_{k-1}, t]$ be given. Let $\{\rho_n\}_{n=1}^\infty$ be a sequence in $(\sigma_{k-1}, \rho)$ satisfying $\rho_n \uparrow \rho$. By definition, $C^\phi(t) \geq (\phi(\rho_n) - a)^+ \wedge \inf_{u \in [\rho_n, t]} \phi(u)$, and letting $n \to \infty$, we obtain

$$(3.9) \quad C^\phi(t) \geq (\phi(\rho-) - a)^+ \wedge \phi(\rho-) \wedge \inf_{u \in [\rho, t]} \phi(u), \qquad \sigma_{k-1} < \rho \leq t < \tau_k.$$

Now let $t \in [\sigma_{k-1}, \tau_k)$ be given. Then there exists $\rho_t$ such that either

$$(3.10) \qquad \rho_t \in [\sigma_{k-1}, t] \quad \text{and} \quad \sup_{s \in [\sigma_{k-1}, t]} \phi(s) = \phi(\rho_t),$$

or else

$$(3.11) \qquad \rho_t \in (\sigma_{k-1}, t] \quad \text{and} \quad \sup_{s \in [\sigma_{k-1}, t]} \phi(s) = \phi(\rho_t-).$$

If (3.11) is the case, which can happen only if $t > \sigma_{k-1}$, then for $u \in [\rho_t, t]$, $\sup_{s \in [\sigma_{k-1}, u]} \phi(s) = \phi(\rho_t-)$ and so (3.4) implies

$$\phi(\rho_t-) - a = (\phi(\rho_t-) - a) \wedge \sup_{s \in [\sigma_{k-1}, u]}(\phi(s) - a) \leq \phi(\rho_t-) \wedge \phi(u),$$

which yields $\phi(\rho_t-) - a \leq \phi(\rho_t-) \wedge \inf_{u \in [\rho_t, t]} \phi(u)$. This inequality together with (3.9) and (3.11) shows that

$$(3.12) \qquad C^\phi(t) \geq (\phi(\rho_t-) - a)^+ = \sup_{s \in [\sigma_{k-1}, t]}(\phi(s) - a)^+.$$

If, on the other hand, (3.10) is the case, then (3.4) implies

$$\phi(\rho_t) - a = \sup_{s \in [\sigma_{k-1}, u]} \phi(s) - a < \phi(u) \qquad \forall u \in [\rho_t, t],$$

and hence $\phi(\rho_t) - a \leq \inf_{u \in [\rho_t, t]} \phi(u)$. This shows that

$$C^\phi(t) \geq (\phi(\rho_t) - a)^+ \wedge \inf_{u \in [\rho_t, t]} \phi(u) = (\phi(\rho_t) - a)^+ = \sup_{s \in [\sigma_{k-1}, t]}(\phi(s) - a)^+.$$

We again have the lower bound (3.12).

To obtain the reverse of inequality (3.12), we consider separately the cases $k = 1$ and $k \geq 2$. If $k = 1$, then $(\phi(s) - a)^+ = 0$ for $s \in [0, \sigma_0)$ and for $t \in (\sigma_0, \tau_1)$,

$$C^\phi(t) \doteq \sup_{s \in [0,t]} \left[ (\phi(s) - a)^+ \wedge \inf_{u \in [s,t]} \phi(u) \right] \leq \sup_{s \in [\sigma_0,t]} (\phi(s) - a)^+,$$

as desired. If $k \geq 2$, we may write $C^\phi(t) = S_1 \vee S_2 \vee S_3$, where

$$(3.13) \qquad S_1 = \sup_{s \in [0, \tau_{k-1}]} \left[ (\phi(s) - a)^+ \wedge \inf_{u \in [s,t]} \phi(u) \right],$$

$$(3.14) \qquad S_2 = \sup_{s \in (\tau_{k-1}, \sigma_{k-1})} \left[ (\phi(s) - a)^+ \wedge \inf_{u \in [s,t]} \phi(u) \right],$$

$$(3.15) \qquad S_3 = \sup_{s \in [\sigma_{k-1}, t]} \left[ (\phi(s) - a)^+ \wedge \inf_{u \in [s,t]} \phi(u) \right].$$

We show that each of the terms $S_i$ is dominated by $\sup_{s \in [\sigma_{k-1},t]} (\phi(s) - a)^+$. For $S_3$, this is obvious. For $S_1$, we use (3.7) and the fact that $t \geq \sigma_{k-1}$ to write

$$(3.16) \qquad \begin{aligned} S_1 &\leq \sup_{s \in [0, \tau_{k-1}]} \inf_{u \in [s,t]} \phi(u) \leq \inf_{u \in [\tau_{k-1}, \sigma_{k-1}]} \phi(u) \\ &\leq \phi(\sigma_{k-1}) - a \leq \sup_{s \in [\sigma_{k-1}, t]} (\phi(s) - a)^+. \end{aligned}$$

Finally, for $s \in (\tau_{k-1}, \sigma_{k-1})$, (3.6) implies $\phi(s) - a < \inf_{u \in [\tau_{k-1}, s]} \phi(u)$, and hence

$$\begin{aligned} S_2 &\leq \sup_{s \in (\tau_{k-1}, \sigma_{k-1})} \left[ \inf_{u \in [\tau_{k-1}, s]} \phi(u) \wedge \inf_{u \in [s,t]} \phi(u) \right] \\ &= \inf_{u \in [\tau_{k-1}, t]} \phi(u) \\ &\leq \inf_{u \in [\tau_{k-1}, \sigma_{k-1}]} \phi(u). \end{aligned}$$

We conclude as in (3.16). □

PROPOSITION 3.3. *We have $C^\phi(t) = 0$ for $t \in [0, \sigma_0)$. For $k \geq 1$, $C^\phi(t) = \inf_{u \in [\tau_k, t]} \phi(u)$ for all $t \in [\tau_k, \sigma_k)$.*

PROOF. Since $\phi \geq 0$, it follows immediately from (2.5) that $C^\phi(t) = 0$ for $t \in [0, \sigma_0)$. Now let $k \geq 1$ and $t \in [\tau_k, \sigma_k)$ be given. By definition,

$$C^\phi(t) = \sup_{s \in [0, \tau_k]} \left[ (\phi(s) - a)^+ \wedge \inf_{u \in [s,t]} \phi(u) \right]$$

(3.17)
$$\vee \sup_{s\in[\tau_k,t]}\left[(\phi(s)-a)^+ \wedge \inf_{u\in[s,t]}\phi(u)\right].$$

It is obvious that
$$\sup_{s\in[0,\tau_k]}\left[(\phi(s)-a)^+ \wedge \inf_{u\in[s,t]}\phi(u)\right] \leq \sup_{s\in[0,\tau_k]}\inf_{u\in[s,t]}\phi(u) = \inf_{u\in[\tau_k,t]}\phi(u).$$

In addition, (3.6) implies
$$\sup_{s\in[\tau_k,t]}\left[(\phi(s)-a)^+ \wedge \inf_{u\in[s,t]}\phi(u)\right]$$
$$\leq \sup_{s\in[\tau_k,t]}\left[\inf_{u\in[\tau_k,s]}\phi(u) \wedge \inf_{u\in[s,t]}\phi(u)\right]$$
$$= \inf_{u\in[\tau_k,t]}\phi(u).$$

We have obtained the upper bound
(3.18) $$C^\phi(t) \leq \inf_{u\in[\tau_k,t]}\phi(u).$$

For the reverse inequality, we observe that there exists $\rho$ such that either

(3.19) $$\rho \in [\sigma_{k-1},\tau_k] \quad\text{and}\quad \sup_{s\in[\sigma_{k-1},\tau_k]}\phi(s) = \phi(\rho),$$

or else

(3.20) $$\rho \in (\sigma_{k-1},\tau_k] \quad\text{and}\quad \sup_{s\in[\sigma_{k-1},\tau_k]}\phi(s) = \phi(\rho-).$$

In either case, we have from (3.4) that for $u \in [\rho,\tau_k)$,
$$\phi(u) > \sup_{s\in[\sigma_{k-1},u]}\phi(s) - a = \sup_{s\in[\sigma_{k-1},\tau_k]}\phi(s) - a,$$

and hence, by (3.5),
(3.21) $$\inf_{u\in[\rho,\tau_k)}\phi(u) \geq \sup_{s\in[\sigma_{k-1},\tau_k]}\phi(s) - a \geq \phi(\tau_k).$$

In the case (3.19), we write
$$C^\phi(t) \geq (\phi(\rho)-a)^+ \wedge \inf_{u\in[\rho,\tau_k)}\phi(u) \wedge \inf_{u\in[\tau_k,t]}\phi(u)$$

and use (3.19), (3.5), and (3.21) to conclude that
(3.22) $$C^\phi(t) \geq \inf_{u\in[\tau_k,t]}\phi(u).$$

In the case (3.20), we choose a sequence $\{\rho_n\}_{n=1}^\infty$ in $(\sigma_{k-1}, \rho)$ with $\rho_n \uparrow \rho$ and write

$$(3.23) \qquad C^\phi(t) \geq (\phi(\rho_n) - a)^+ \wedge \inf_{u \in [\rho_n, \tau_k)} \phi(u) \wedge \inf_{u \in [\tau_k, t]} \phi(u).$$

Letting $n \to \infty$, we obtain

$$C^\phi(t) \geq (\phi(\rho-) - a)^+ \wedge \phi(\rho-) \wedge \inf_{u \in [\rho, \tau_k)} \phi(u) \wedge \inf_{u \in [\tau_k, t]} \phi(u)$$

$$\geq (\phi(\rho-) - a)^+ \wedge \inf_{u \in [\rho, \tau_k)} \phi(u) \wedge \inf_{u \in [\tau_k, t]} \phi(u).$$

We now use (3.20), (3.5), and (3.21) to conclude (3.22). □

In summary, Propositions 3.2 and 3.3 imply that $C^\phi(t)$ given by (2.5) has the form

$$(3.24) \quad C^\phi(t) = \begin{cases} 0, & \text{if } 0 \leq t < \sigma_0, \\ \sup_{s \in [\sigma_{k-1}, t]} (\phi(s) - a)^+, & \text{if } \sigma_{k-1} \leq t < \tau_k, k \geq 1, \\ \inf_{u \in [\tau_k, t]} \phi(u), & \text{if } \tau_k \leq t < \sigma_k, k \geq 1. \end{cases}$$

The inequalities (3.5) and (3.7) imply $\sup_{s \in [\sigma_{k-1}, \tau_k]} \phi(s) = \sup_{s \in [\sigma_{k-1}, \tau_k)} \phi(s)$ and $\inf_{u \in [\tau_k, \sigma_k]} \phi(u) = \inf_{u \in [\tau_k, \sigma_k)} \phi(u)$. Moreover, when combined with (3.24) and the fact that $\phi \geq 0$, these inequalities show that for $k \geq 1$,

$$(3.25) \qquad C^\phi(\tau_k-) = \sup_{s \in [\sigma_{k-1}, \tau_k)} (\phi(s) - a)^+ \geq \phi(\tau_k) = C^\phi(\tau_k),$$

$$(3.26) \qquad C^\phi(\sigma_k-) = \inf_{u \in [\tau_k, \sigma_k)} \phi(u) \leq \phi(\sigma_k) - a = C^\phi(\sigma_k).$$

We define $C^\phi(0-) = 0$ and we have

$$(3.27) \qquad C^\phi(\sigma_0-) = 0 \leq C^\phi(\sigma_0) = (\phi(\sigma_0) - a)^+ = \phi(\sigma_0) - a,$$

where the last equality holds due to (3.8). In particular, $C^\phi$ is increasing on each interval $[\sigma_{k-1}, \tau_k)$, with a possible upward jump at $\sigma_{k-1}$, and $C^\phi$ is decreasing on each interval $[\tau_k, \sigma_k)$, with a possible downward jump at $\tau_k$.

THEOREM 3.4. *Let $\phi \in \mathcal{D}_+[0, \infty)$ be given, define $C^\phi$ by (2.5), and set $\bar{\phi} = \phi - C^\phi$. Then $C^\phi \in \mathcal{BV}[0, \infty)$, $\bar{\phi} \in \mathcal{D}[0, \infty)$, and $\bar{\phi}$ takes values only in $[0, a]$. Furthermore,*

$$(3.28) \qquad |C^\phi|(t) = \int_0^t \mathbb{I}_{\{\bar{\phi}(s) = 0 \text{ or } \bar{\phi}(s) = a\}} \, d|C^\phi|(s),$$

$$(3.29) \qquad C^\phi(t) = -\int_0^t \mathbb{I}_{\{\bar{\phi}(s) = 0\}} \, d|C^\phi|(s) + \int_0^t \mathbb{I}_{\{\bar{\phi}(s) = a\}} \, d|C^\phi|(s).$$

PROOF. From (3.24) we see that $C^\phi \in \mathcal{BV}[0,\infty)$. From its definition (2.5), we see that $C^\phi$ further satisfies $(\phi - a)^+ \leq C^\phi \leq \phi$, and hence

$$0 \leq \bar{\phi} \leq a \wedge \phi. \tag{3.30}$$

Moreover, the rightmost equalities in the relations (3.25)–(3.27) show that

$$\bar{\phi}(\tau_k) = 0 \quad \text{and} \quad \bar{\phi}(\sigma_{k-1}) = a, \qquad k \geq 1. \tag{3.31}$$

Since $C^\phi = 0$ on $[0, \sigma_0)$, we only need to consider $t \geq \sigma_0$ in what follows. Define the set

$$A \doteq \{t \geq \sigma_0 : \bar{\phi}(t) \in (0, a)\}. \tag{3.32}$$

We show below that $\int_A d|C^\phi| = 0$, so that (3.28) holds. We further show that for $t \geq \sigma_0$,

$$\bar{\phi}(t) = 0 \implies t \in [\tau_k, \sigma_k) \quad \text{for some } k, \tag{3.33}$$

whereas

$$\bar{\phi}(t) = a \implies t \in [\sigma_{k-1}, \tau_k) \quad \text{for some } k. \tag{3.34}$$

We can then conclude that $C^\phi$ does not increase on $\{t \geq 0 | \bar{\phi}(t) = 0\}$ (the positive variation of $C^\phi$ assigns zero measure to this set) and $C^\phi$ does not decrease on the set $\{t \geq 0 | \bar{\phi}(t) = a\}$ (the negative variation of $C^\phi$ assigns zero measure to this set). This together with (3.28) will imply (3.29).

We first establish (3.33) and (3.34). Suppose $t \in [\sigma_{k-1}, \tau_k)$ for some $k$. Then (3.4) and either (3.7) or (3.8) imply

$$\phi(t) > \sup_{s \in [\sigma_{k-1}, t]} \phi(s) - a \geq \phi(\sigma_{k-1}) - a \geq 0.$$

From this and (3.24) we have

$$C^\phi(t) = \sup_{s \in [\sigma_{k-1}, t]} (\phi(s) - a)^+ = \sup_{s \in [\sigma_{k-1}, t]} \phi(s) - a < \phi(t).$$

Therefore, $\bar{\phi}(t) = \phi(t) - C^\phi(t) > 0$. This is the contrapositive of (3.33). Similarly, suppose $t \in [\tau_k, \sigma_k)$ for some $k$. Then (3.24) and (3.6) imply $C^\phi(t) = \inf_{u \in [\tau_k, t]} \phi(u) > \phi(t) - a$, so that $\bar{\phi}(t) = \phi(t) - C^\phi(t) < a$. This is the contrapositive of (3.34).

We next show that $\int_A d|C^\phi| = 0$. For $t \in A$, define

$$\alpha(t) \doteq \inf\{s \in [\sigma_0, t] | (s, t] \subset A\}, \qquad \beta(t) \doteq \sup\{s \in [t, \infty) | [t, s) \in A\}.$$

Because of the right-continuity of $\bar{\phi}$, we have $\beta(t) \notin A$, whereas $\alpha(t)$ might or might not be in $A$. We also have $\alpha(t) \leq t < \beta(t)$, and so the open interval $(\alpha(t), \beta(t))$ is nonempty. It follows that $A$ is the countable union of such disjoint open intervals together with a countable set of left endpoints, that is,

$$A = \left(\bigcup_{i \in I}(\alpha_i, \beta_i)\right) \cup \{\alpha_j | j \in J\},$$

where $I$ is a countable index set and $J \subset I$.

As a first step in showing $\int_A d|C^\phi| = 0$, we show that if $j \in J$, so $\alpha_j \in A$, then $C^\phi$ is continuous at $\alpha_j$. From (3.31) we see that $\alpha_j$ is in the interior of an interval of the form $(\tau_k, \sigma_k)$ or of the form $(\sigma_{k-1}, \tau_k)$. By the definition of $\alpha_j$, there is a sequence of points $\{\gamma_n\}_{n=1}^\infty$ in $(0, \alpha_j) \cap A^c$ such that $\gamma_n \uparrow \alpha_j$.

We consider first the case that $\bar{\phi}(\gamma_n) = a$, or equivalently, $C^\phi(\gamma_n) = \phi(\gamma_n) - a$, for infinitely many values of $n$. From (3.34), we see that $\gamma_n \in [\sigma_{k-1}, \tau_k)$ for some $k$. By choosing $n$ sufficiently large, we may assume that $k$ does not depend on $n$ and $\alpha_j \in (\sigma_{k-1}, \tau_k)$. We have

$$a = \phi(\gamma_n) - C^\phi(\gamma_n) = \phi(\gamma_n) - \sup_{s \in [\sigma_{k-1}, \gamma_n]}(\phi(s) - a)^+$$

$$\leq \phi(\gamma_n) - (\phi(\gamma_n) - a)^+ = \phi(\gamma_n) \wedge a$$

$$\leq a.$$

Therefore, the above inequalities must be equalities and we conclude that

$$0 \leq \phi(\gamma_n) - a = C^\phi(\gamma_n) = \sup_{s \in [\sigma_{k-1}, \gamma_n]}(\phi(s) - a)^+.$$

Letting $n \to \infty$, we see that

$$0 \leq \phi(\alpha_j-) - a = C^\phi(\alpha_j-) = \sup_{s \in [\sigma_{k-1}, \alpha_j)}(\phi(s) - a)^+.$$

On the other hand, $C^\phi(\alpha_j) = \sup_{s \in [\sigma_{k-1}, \alpha_j]}(\phi(s) - a)^+$. This shows that $C^\phi(\alpha_j) \geq C^\phi(\alpha_j-)$. Furthermore, $C^\phi(\alpha_j) > C^\phi(\alpha_j-)$ implies $C^\phi(\alpha_j) = \phi(\alpha_j) - a$. But in this case, $\bar{\phi}(\alpha_j) = a$. This contradicts the membership of $\alpha_j$ in $A$ and establishes the continuity of $C^\phi$ at $\alpha_j$.

If $\bar{\phi}(\gamma_n) = a$ does not hold for infinitely many values of $n$, then $\bar{\phi}(\gamma_n) = 0$, or equivalently, $C^\phi(\gamma_n) = \phi(\gamma_n)$, must hold for infinitely many values of $n$. From (3.33), we see that $\gamma_n \in [\tau_k, \sigma_k)$ for some $k$. By choosing $n$ sufficiently large, we may assume that $k$ does not depend on $n$ and $\alpha_j \in (\tau_k, \sigma_k)$. We have

$$0 = \phi(\gamma_n) - C^\phi(\gamma_n) = \phi(\gamma_n) - \inf_{u \in [\tau_k, \gamma_n]}\phi(u) \geq 0.$$

Therefore, the above inequality must be an equality and we conclude that

$$\phi(\gamma_n) = C^\phi(\gamma_n) = \inf_{u \in [\tau_k, \gamma_n]}\phi(u).$$

Letting $n \to \infty$, we see that
$$\phi(\alpha_j-) = C^\phi(\alpha_j-) = \inf_{u \in [\tau_k, \alpha_j)} \phi(u).$$

On the other hand, $C^\phi(\alpha_j) = \inf_{u \in [\tau_k, \alpha_j]} \phi(u)$. This shows that $C^\phi(\alpha_j) \leq C^\phi(\alpha_j-)$. Furthermore, $C^\phi(\alpha_j) < C^\phi(\alpha_j-)$ implies $C^\phi(\alpha_j) = \phi(\alpha_j)$. But in this case, $\bar\phi(\alpha_j) = 0$. This contradicts the membership of $\alpha_j$ in $A$, which establishes the continuity of $C^\phi$ at $\alpha_j$.

To establish $\int_A d|C^\phi| = 0$, it remains only to show that $\int_{(\alpha_i,\beta_i)} d|C^\phi| = 0$ for every $i \in I$. Because $\bar\phi$ is strictly between 0 and $a$ on $(\alpha_i, \beta_i)$, (3.31) shows that $(\alpha_i, \beta_i)$ must be entirely contained in an interval of the form $(\tau_k, \sigma_k)$ or of the form $(\sigma_{k-1}, \tau_k)$. We consider the latter case; the former case is analogous. It suffices to show that $C^\phi$ is constant on $[a_i, b_i]$ whenever $\alpha_i < a_i < b_i < \beta_i$, where
$$C^\phi(t) = \sup_{s \in [\sigma_{k-1}, t]} (\phi(s) - a)^+ \qquad \forall t \in (\alpha_i, \beta_i).$$

Define
$$\rho = \inf\{t \in [a_i, b_i] | C^\phi(t) > C^\phi(a_i)\}.$$

Assume $\rho < \infty$. Because $C^\phi$ is right-continuous, we must have $C^\phi(t) = C^\phi(a_i)$ for all $t \in [a_i, \rho)$ and either $C^\phi(\rho) = \phi(\rho) - a > C(a_i)$ or else $C^\phi(\rho) = \phi(\rho) - a = C^\phi(a_i)$. In either case, $\bar\phi(\rho) = a$, contradicting the definition of $A$. Therefore, $\rho = \infty$ and $C^\phi$ is constant on $[a_i, b_i]$. □

PROOF OF THEOREM 1.4. Let $\psi \in \mathcal{D}[0,\infty)$ be given and define $\phi = \Gamma_0(\psi)$. Then $\eta \doteq \phi - \psi \in \mathcal{I}[0,\infty)$ satisfies [see (1.2)]

$$(3.35) \quad \eta(t) = \int_0^t \mathbb{I}_{\{\phi(s)=0\}} \, d\eta(s), \qquad \int_0^t \mathbb{I}_{\{\phi(s)>0\}} \, d\eta(s) = 0 \qquad \forall t \geq 0.$$

With $C^\phi$ defined by (2.5), set
$$\bar\phi = \Lambda_a(\phi) = \phi - C^\phi = \psi + \eta - C^\phi.$$

Theorem 3.4 implies $\eta - C^\phi \in \mathcal{BV}[0,\infty)$, $\bar\phi \in \mathcal{D}[0,\infty)$, and $\bar\phi$ takes values only in $[0,a]$. It remains to show that for all $t \geq 0$,

$$(3.36) \quad |\eta - C^\phi|(t) = \int_0^t \mathbb{I}_{\{\bar\phi(s)=0 \text{ or } \bar\phi(s)=a\}} \, d|\eta - C^\phi|(s),$$

$$(3.37) \quad \eta(t) - C^\phi(t) = \int_0^t \mathbb{I}_{\{\bar\phi(s)=0\}} \, d|\eta - C^\phi|(s) - \int_0^t \mathbb{I}_{\{\bar\phi(s)=a\}} \, d|\eta - C^\phi|(s).$$

Because $\{s|\phi(s) = 0\} \subset \{s|\bar\phi(s) = 0\}$ [see (3.30)] and $C^\phi$ is decreasing on this set [see (3.29)], (3.35) implies $|\eta - C^\phi| = \eta + |C^\phi|$. Equations (3.36) and (3.37) follow from (3.35), (3.28) and (3.29). □

We now present the proof of Corollary 1.6.

PROOF OF COROLLARY 1.6. We first prove (1.15) for $i = \infty$. For $\phi_1, \phi_2 \in \mathcal{D}[0, T]$, we have

$$(3.38) \qquad \|\Lambda_a(\phi_1) - \Lambda_a(\phi_2)\|_T \leq \|\phi_1 - \phi_2\|_T + \|C^{\phi_1} - C^{\phi_2}\|_T.$$

For $t \in [0, T]$, because $(a_1 \wedge b_1) - (a_2 \wedge b_2) \leq (a_1 - a_2) \vee (b_1 - b_2)$, we have

$$\begin{aligned}
C^{\phi_1}(t) &- C^{\phi_2}(t) \\
&\leq \sup_{s \in [0, t]} [R_t(\phi_1)(s) - R_t(\phi_2)(s)] \\
&\leq \sup_{s \in [0, t]} \left[ |(\phi_1(s) - a)^+ - (\phi_2(s) - a)^+| \vee \left| \inf_{u \in [s, t]} \phi_1(u) - \inf_{u \in [s, t]} \phi_2(u) \right| \right] \\
&\leq \sup_{s \in [0, t]} \left[ |\phi_1(s) - \phi_2(s)| \vee \sup_{u \in [s, t]} |\phi_1(u) - \phi_2(u)| \right] \\
&\leq \|\phi_1 - \phi_2\|_T.
\end{aligned}$$

Taking the supremum over $t \in [0, T]$ and interchanging $\phi_1$ and $\phi_2$, we get

$$(3.39) \qquad \|C^{\phi_1} - C^{\phi_2}\|_T \leq \|\phi_1 - \phi_2\|_T.$$

From (3.38) and (3.39), we obtain (1.15) for $i = \infty$.

Now let $\mathcal{M}$ be the class of strictly increasing continuous functions $\lambda$ of $[0, T]$ onto itself. Then for any $\lambda \in \mathcal{M}$, the scaling property

$$(3.40) \qquad \Lambda_a(\phi \circ \lambda) = \Lambda_a(\phi) \circ \lambda$$

is easily deduced directly from the definition of $\Lambda_a$. Moreover, by the definition of $d_0$, given any $\phi_1, \phi_2 \in \mathcal{D}[0, T]$, $\phi_1 \neq \phi_2$, for every $\delta > 0$ there exists $\lambda \in \mathcal{M}$ (possibly depending on $\delta$) such that

$$\sup_{t \in [0, T]} |\lambda(t) - t| \leq d_0(\phi_1, \phi_2) + \delta[1 \wedge d_0(\phi_1, \phi_2)]$$

and

$$\sup_{t \in [0, T]} |\phi_1(t) - \phi_2(\lambda(t))| \leq d_0(\phi_1, \phi_2) + \delta[1 \wedge d_0(\phi_1, \phi_2)].$$

The scaling property (3.40) along with (1.15) for $i = \infty$ implies that

$$\sup_{t \in [0, T]} |\Lambda_a(\phi_1)(t) - \Lambda_a(\phi_2)(\lambda(t))| \leq 2(d_0(\phi_1, \phi_2) + \delta[1 \wedge d_0(\phi_1, \phi_2)]).$$

Since this is true for all $\delta > 0$, by the definition of $d_0$ this implies that

$$d_0(\Lambda_a(\phi_1), \Lambda_a(\phi_2)) \leq 2 d_0(\phi_1, \phi_2),$$

which is the inequality (1.15) for $i = 0$. Clearly, (1.15) holds also in the case $i = 0$, $\phi_1 = \phi_2 \in \mathcal{D}[0, T]$.

We now prove (1.15) for $i = 1$. For a given $\phi \in \mathcal{D}[0, T]$, let $\phi(0-) \doteq \phi(0)$ and let

$$(3.41) \quad G_\phi = \{(t, z) \in [0, T] \times \mathbb{R} : z \in [\phi(t-) \wedge \phi(t), \phi(t-) \vee \phi(t)]\}$$

be the graph of $\phi$ ordered by the following relation: $(t_1, z_1) \leq (t_2, z_2)$ if either $t_1 < t_2$ or $t_1 = t_2$ and $|\phi(t_1-) - z_1| \leq |\phi(t_1-) - z_2|$. Let $\Pi(\phi)$ be the set of all parametric representations of $G_\phi$, that is, continuous nondecreasing (in the order relation just defined) functions $(r, g)$ mapping $[0, 1]$ onto $G_\phi$. For $\phi_1, \phi_2 \in \mathcal{D}[0, T]$,

$$d_1(\phi_1, \phi_2) \doteq \inf\{\|r_1 - r_2\|_T \vee \|g_1 - g_2\|_T : (r_i, g_i) \in \Pi(\phi_i), i = 1, 2\}.$$

We show in Lemma A.1 in the Appendix that if $(r, g) \in \Pi(\phi)$, then $(r, \Lambda_a(g)) \in \Pi(\Lambda_a(\phi))$. Therefore,

$$\begin{aligned} d_1(\Lambda_a(\phi_1), \Lambda_a(\phi_2)) \\ &\leq \inf\{\|r_1 - r_2\|_T \vee \|\Lambda_a(g_1) - \Lambda_a(g_2)\|_T : (r_i, g_i) \in \Pi(\phi_i), i = 1, 2\} \\ &\leq 2d_1(\phi_1, \phi_2), \end{aligned}$$

where the last inequality follows from (1.15) for $i = \infty$. We have proved (1.15).

It is well known (see, e.g., Lemma 13.5.1 and Theorem 13.5.1 of [16]) that for any $T < \infty$ and $\psi_1, \psi_2 \in \mathcal{D}[0, T]$

$$(3.42) \qquad d_i(\Gamma_0(\psi_1), \Gamma_0(\psi_2)) \leq 2d_i(\psi_1, \psi_2)$$

for $i = 0, 1, \infty$. The representation $\Gamma_{0,a} = \Lambda_a \circ \Gamma_0$ stated in (1.13), along with (1.15) and (3.42) then implies that (1.16) holds with $L = 4$.

By the argument in Theorem 12.9.4 in [16], the validity of (1.15) and (1.16) on $\mathcal{D}[0, T]$ for every $T > 0$ implies the same bound on $\mathcal{D}[0, \infty)$. $\square$

REMARKS. Example 13.5.1 in [16] shows that the bound in (3.42) with $i = \infty$ is tight. Similarly, the bound (1.15) for $i = \infty$ is tight. To see this, let us consider $\phi_1, \phi_2 \in \mathcal{D}[0, 1]$ defined by $\phi_1 = 2\mathbb{I}_{[0,1]}$, $\phi_2 = 3\mathbb{I}_{[0,1/2)} + \mathbb{I}_{[1/2,1]}$. With $a = 2$ we have $\Lambda_a(\phi_1) = \phi_1$, $\Lambda_a(\phi_2) = 2\mathbb{I}_{[0,1/2)}$, $\|\phi_1 - \phi_2\|_1 = 1$ and $\|\Lambda_a(\phi_1) - \Lambda_a(\phi_2)\|_1 = 2$. However, Theorem 14.8.1 in [16] shows that (1.16) for $i = \infty$ (and thus also for $i = 0, 1$) actually holds with $L = 2$. Clearly, the bound (1.16) with $L = 2$ is tight, because the bound (3.42) is tight.

**4. Comparison properties of the double reflection map.** In this section we present the proof of Theorem 1.7. We first establish some preliminary results that may be of independent interest. In the proofs we make repeated use of the elementary inequalities $[b_1 + b_2]^+ \leq b_1^+ + b_2^+$ and $[b_1 - b_2]^+ \geq b_1^+ - b_2^+$ for $b_1, b_2 \in \mathbb{R}$, without explicit reference.

LEMMA 4.1. *Given $c_0, c_0' \in \mathbb{R}$ and $\psi, \psi' \in \mathcal{D}[0, \infty)$ with $\psi(0) = \psi'(0) = 0$, suppose $(\phi, \eta)$ and $(\phi', \eta')$ solve the Skorokhod problem on $[0, \infty)$ for $c_0 + \psi$ and $c_0' + \psi'$, respectively. If there exists $\nu \in \mathcal{I}[0, \infty)$ such that $\psi' \leq \psi \leq \psi' + \nu$, then the following two properties are satisfied:*

1. $\eta - [c_0' - c_0]^+ \leq \eta' \leq \eta + \nu + [c_0 - c_0']^+$;
2. $\phi' - \nu - [c_0' - c_0]^+ \leq \phi \leq \phi' + \nu + [c_0 - c_0']^+$.

*Moreover, if $\psi = \psi' + \nu$ then*

(4.1) $$\phi' - [c_0' - c_0]^+ \leq \phi \leq \phi' + \nu + [c_0 - c_0']^+.$$

PROOF. Using the explicit representations for $\eta$ and $\eta'$ that follow from (1.1), along with the fact that $\nu \in \mathcal{I}[0, \infty)$ and $\psi \leq \psi' + \nu$, we see that for every $t \in [0, \infty)$,

$$\begin{aligned}
\eta(t) &= \sup_{s \in [0,t]} [-c_0 - \psi(s)]^+ \\
&\geq \sup_{s \in [0,t]} [-c_0' - \psi'(s) - \nu(s) - c_0 + c_0']^+ \\
&\geq \sup_{s \in [0,t]} [-c_0' - \psi'(s) - \nu(t) - c_0 + c_0']^+ \\
&\geq \sup_{s \in [0,t]} [-c_0' - \psi'(s)]^+ - [\nu(t) + c_0 - c_0']^+ \\
&\geq \eta'(t) - \nu(t) - [c_0 - c_0']^+.
\end{aligned}$$

Likewise, (1.1) and the fact that $\psi \geq \psi'$ shows that for every $t \in [0, \infty)$,

(4.2) $$\begin{aligned}
\eta'(t) &= \sup_{s \in [0,t]} [-c_0' - \psi'(s)]^+ \\
&\geq \sup_{s \in [0,t]} [-c_0 - \psi(s) - (c_0' - c_0)]^+ \\
&\geq \sup_{s \in [0,t]} [-c_0 - \psi(s)]^+ - [c_0' - c_0]^+ \\
&= \eta(t) - [c_0' - c_0]^+.
\end{aligned}$$

When combined, the last two relations establish property 1. Moreover, the first relation and the fact that $\eta' = -c_0' - \psi' + \phi'$ also implies that

$$\phi = \psi + c_0 + \eta \geq \psi + c_0 - c_0' - \psi' + \phi' - \nu - [c_0 - c_0']^+$$

$$= \phi' + \psi - \psi' - \nu - [c'_0 - c_0]^+,$$

which is no less than $\phi' - \nu - [c'_0 - c_0]^+$ if $\psi' \leq \psi \leq \psi' + \nu$ and is no less than $\phi' - [c'_0 - c_0]^+$ if $\psi = \psi' + \nu$. On the other hand, the second relation, (4.2), shows that

$$\phi = c_0 + \psi + \eta \leq c'_0 + \psi' + \eta' + c_0 - c'_0 + [c'_0 - c_0]^+ + \psi - \psi'$$
$$= \phi' + [c_0 - c'_0]^+ + \psi - \psi'.$$

Together, the last two displays establish property 2 and (4.1). □

The representation (1.13) for $\Gamma_{0,a}$ as the composition of $\Lambda_a$ and $\Gamma_0$, allows us to easily deduce the following corollary from Lemma 4.1.

COROLLARY 4.2. *Given $a > 0$, $c_0, c'_0 \in \mathbb{R}$ and $\psi, \psi' \in \mathcal{D}[0, \infty)$ with $\psi(0) = \psi'(0) = 0$, suppose $(\bar{\phi}, \bar{\eta})$ and $(\bar{\phi}', \bar{\eta}')$ solve the Skorokhod problem on $[0, a]$ for $c_0 + \psi$ and $c'_0 + \psi'$, respectively. If $\psi = \psi' + \nu$, where $\nu \in \mathcal{I}[0, \infty)$, then the following two properties hold:*

1. $\bar{\eta} - 2[c'_0 - c_0]^+ \leq \bar{\eta}' \leq \bar{\eta} + 2\nu + 2[c_0 - c'_0]^+$;
2. $[-|c'_0 - c_0| - \nu] \vee [-a] \leq \bar{\phi}' - \bar{\phi} \leq [|c'_0 - c_0| + \nu] \wedge a.$

PROOF. Let $C = C^\phi$ be the function defined in (2.5) and let $C' = C^{\phi'}$. From the first inequality in (4.1) of Lemma 4.1, it follows that

$$C'(t) = \sup_{s \in [0,t]} \left[ (\phi'(s) - a)^+ \wedge \inf_{u \in [s,t]} \phi'(u) \right]$$

$$\leq \sup_{s \in [0,t]} \left[ (\phi(s) - a + [c'_0 - c_0]^+)^+ \wedge \inf_{u \in [s,t]} (\phi(u) + [c'_0 - c_0]^+) \right]$$

$$\leq \sup_{s \in [0,t]} \left[ (\phi(s) - a)^+ \wedge \inf_{u \in [s,t]} \phi(u) \right] + [c'_0 - c_0]^+ = C(t) + [c'_0 - c_0]^+.$$

Similarly, the second inequality in (4.1) along with the fact that $\nu$ is nondecreasing implies that $C'(t)$ is equal to

$$\sup_{s \in [0,t]} \left[ (\phi'(s) - a)^+ \wedge \inf_{u \in [s,t]} \phi'(u) \right]$$

$$\geq \sup_{s \in [0,t]} \left[ (\phi(s) - a - \nu(t) - [c_0 - c'_0]^+)^+ \right.$$
$$\left. \wedge \inf_{u \in [s,t]} (\phi(u) - \nu(t) - [c_0 - c'_0]^+) \right]$$

$$\geq \sup_{s \in [0,t]} \left[ (\phi(s) - a)^+ \wedge \inf_{u \in [s,t]} \phi(u) \right] - \nu(t) - [c_0 - c'_0]^+$$

$$= C(t) - \nu(t) - [c_0 - c'_0]^+.$$

Let $\eta = \Gamma_0(c_0 + \psi) - c_0 - \psi$ and, likewise, let $\eta' = \Gamma_0(c'_0 + \psi') - c'_0 - \psi'$, and note that due to the representation for $\Gamma_{0,a}$ in (1.13), the definition (1.11) of $\Lambda_a$ and the definitions of $C, C'$, we can write $\bar{\eta} = \eta - C$ and $\bar{\eta}' = \eta' - C'$. The last two displays, together with property 1 of Lemma 4.1, then show that

$$\bar{\eta} = \eta - C \leq \eta' + [c'_0 - c_0]^+ - C' + [c'_0 - c_0]^+ = \bar{\eta}' + 2[c'_0 - c_0]^+$$

and

$$\bar{\eta} = \eta - C \geq \eta' - \nu - [c_0 - c'_0]^+ - C' - \nu - [c_0 - c'_0]^+$$
$$= \bar{\eta}' - 2\nu - 2[c_0 - c'_0]^+,$$

which establishes the first property of the corollary. The second property follows from the first property, the fact that $\bar{\phi}', \bar{\phi} \in [0, a]$ and the relation

(4.3) $\quad \bar{\phi}' - \bar{\phi} = c'_0 + \psi' + \bar{\eta}' - c_0 - \psi - \bar{\eta} = c'_0 - c_0 - \nu + \bar{\eta}' - \bar{\eta}.$ $\quad\square$

We introduce the family of shift operators $T_r : \mathcal{D}[0, \infty) \to \mathcal{D}[0, \infty)$, $r \in [0, \infty)$, defined by

$$[T_r f](t) = f(r + t) - f(r) \qquad \text{for } t \in [0, \infty).$$

We shall also make use of the well known (and easily verified) fact that if $\phi = \Gamma(\psi)$, where $\Gamma$ is either the one-sided reflection map at zero or $a$, or the two-sided reflection map on $[0, a]$, then for every $\alpha > 0$,

(4.4) $\qquad\qquad \phi(\alpha + s) = \Gamma(\phi(\alpha) + T_\alpha \psi)(s).$

REMARK 4.3. The first and second inequalities in Corollary 4.2 can be strengthened to the inequalities

(4.5) $\qquad\qquad -[c'_0 - c_0]^+ \leq \bar{\eta}' - \bar{\eta} \leq [c_0 - c'_0]^+ + \nu$

and

(4.6) $\qquad [-[c_0 - c'_0]^+ - \nu] \vee [-a] \leq \bar{\phi}' - \bar{\phi} \leq [c'_0 - c_0]^+ \wedge a,$

which are both easily seen to be tight. Since $\bar{\phi}(t), \bar{\phi}'(t) \in [0, a]$ for all $t \in [0, \infty)$, in order to establish (4.6), it suffices to show that

(4.7) $\quad -[c_0 - c'_0]^+ - \nu(t) \leq \bar{\phi}'(t) - \bar{\phi}(t) \leq [c'_0 - c_0]^+ \qquad \text{for } t \in [0, \infty).$

In order to establish this relation, we use the projection operator $\pi$ of (1.8), which is clearly monotone and Lipschitz with Lipschitz constant 1.

First suppose $c_0 \geq c'_0$. Then, due to the monotonicity property of the projection operator $\pi$ and the Lipschitz continuity of $\Gamma_{0,a}$, Lemma 4.2 of [12] shows that the upper bound $\bar{\phi}' - \bar{\phi} \leq 0 = [c'_0 - c_0]^+$ in (4.7) holds, while the lower bound in (4.7) follows from the first inequality in part 2 of Corollary 4.2.

Now suppose $c_0 < c_0'$. Define

$$\tau \doteq \inf\{t \geq 0 : \bar{\phi}(t) \geq \bar{\phi}'(t)\}.$$

The fact that $\bar{\phi}(0) = \pi(c_0) \leq \pi(c_0') = \bar{\phi}'(0)$ and $\bar{\phi}(t), \bar{\phi}'(t) \in [0, a]$ imply $\bar{\phi}(t) < a$ and $\bar{\phi}'(t) > 0$ for $t \in [0, \tau)$. (It could happen that $\pi(c_0) = \pi(c_0')$, and then $\tau = 0$ and all assertions concerning $t \in [0, \tau)$ are vacuously true.) Definitions 1.1, 1.2 and relation (1.4) then show that for $t \in [0, \tau)$, $\bar{\phi}(t) = \Gamma_0(c_0 + \psi)(t)$ and $\bar{\phi}'(t) = \Gamma_a(c_0' + \psi')(t)$. Therefore for $t \in [0, \tau)$, $c_0 + \psi(t) \leq \bar{\phi}(t) < \bar{\phi}'(t) \leq c_0' + \psi'(t)$, which in turn implies that

$$-\nu(t) \leq 0 \leq \bar{\phi}'(t) - \bar{\phi}(t) \leq c_0' - c_0 + \psi'(t) - \psi(t) \leq c_0' - c_0$$
$$\text{for } t \in [0, \tau).$$

This shows that (4.7) is satisfied for $t \in [0, \tau)$. In particular, this implies that $\bar{\phi}'(\tau-) \geq \bar{\phi}(\tau-) - \nu(\tau-)$. By the monotonicity property of the projection operator $\pi$, we have

$$\begin{aligned}
\bar{\phi}'(\tau) &= \pi(\bar{\phi}'(\tau-) + \psi'(\tau) - \psi'(\tau-)) \\
&\geq \pi(\bar{\phi}(\tau-) - \nu(\tau-) + \psi(\tau) - \psi(\tau-) - (\nu(\tau) - \nu(\tau-))) \\
&\geq \pi(\bar{\phi}(\tau-) + \psi(\tau) - \psi(\tau-)) - \nu(\tau) \\
&= \bar{\phi}(\tau) - \nu(\tau),
\end{aligned}$$
(4.8)

where the explicit definition of $\pi$ is used to obtain the second inequality. Now for $s \in [0, \infty)$, $\bar{\phi}(\tau + s) = \Gamma_{0,a}(\bar{\phi}(\tau) + T_\tau \psi)(s)$ and, likewise, $\bar{\phi}'(\tau + s) = \Gamma_{0,a}(\bar{\phi}'(\tau) + T_\tau \psi')(s)$. Since $\bar{\phi}(\tau) \geq \bar{\phi}'(\tau)$ due to the right-continuity of $\bar{\phi}, \bar{\phi}'$, we can apply (4.7) [with $c_0, c_0', \psi, \psi'$ and $\nu$ replaced by $\bar{\phi}(\tau), \bar{\phi}'(\tau), T_\tau \psi, T_\tau \psi'$ and $T_\tau \nu$], and use (4.8) to obtain for $s \in [0, \infty)$,

$$\begin{aligned}
-\nu(\tau + s) &\leq -[\bar{\phi}(\tau) - \bar{\phi}'(\tau)]^+ - T_\tau \nu(s) \\
&\leq \bar{\phi}'(\tau + s) - \bar{\phi}(\tau + s) \\
&\leq [\bar{\phi}'(\tau) - \bar{\phi}(\tau)]^+ \\
&= 0,
\end{aligned}$$

which shows that (4.7) also holds for $t \in [\tau, \infty)$.

We have established (4.7), and hence (4.6). The inequality (4.5) can be deduced from (4.6) using the basic relation

$$\bar{\eta}' - \bar{\eta} = \bar{\phi}' - \bar{\phi} - (c_0' - c_0) - (\psi' - \psi) = \bar{\phi}' - \bar{\phi} - (c_0' - c_0) + \nu.$$

Although Corollary 4.2 provides bounds on the difference between the net constraining terms $\bar{\eta}$ and $\bar{\eta}'$, it is often desirable to compare the individual constraining terms at the upper and lower barriers. Such bounds are provided in Theorem 1.7. To establish these bounds, we recall that if $(\bar{\phi}, \bar{\eta})$ solves the Skorokhod problem on $[0, a]$ for $\psi \in \mathcal{D}[0, \infty)$, and if $\bar{\eta}$ admits the decomposition $\bar{\eta} = \bar{\eta}_\ell - \bar{\eta}_u$ that satisfies (1.9), then for any $t \in [0, \infty)$,

$$\begin{aligned}
\bar{\eta}_\ell(t) - \bar{\eta}_\ell(t-) &= \sup_{s \in [0,t]} [\bar{\eta}_u(s) - \psi(s)]^+ - \sup_{s \in [0,t)} [\bar{\eta}_u(s) - \psi(s)]^+ \\
(4.9) \qquad &= [\bar{\eta}_u(t) - \psi(t) - \bar{\eta}_\ell(t-)]^+ \\
&= [-\bar{\phi}(t-) - \psi(t) + \psi(t-) + \bar{\eta}_u(t) - \bar{\eta}_u(t-)]^+.
\end{aligned}$$

PROOF OF THEOREM 1.7. Define

$$\alpha \doteq \inf\{t > 0 : \bar{\eta}_\ell(t) + \nu(t) + [c_0 - c'_0]^+ < \bar{\eta}'_\ell(t) \text{ or } \bar{\eta}_u(t) + [c'_0 - c_0]^+ < \bar{\eta}'_u(t)\},$$

with $\alpha \doteq \infty$ if the infimum is over the empty set. Then the definition of $\alpha$ dictates that the following two relations are satisfied for $s \in [0, \alpha)$:

$$(4.10) \qquad \bar{\eta}'_\ell(s) \le \bar{\eta}_\ell(s) + \nu(s) + [c_0 - c'_0]^+;$$

$$(4.11) \qquad \bar{\eta}'_u(s) \le \bar{\eta}_u(s) + [c'_0 - c_0]^+.$$

Suppose $\alpha < \infty$. Then we claim (and prove below) that it is also true that

$$(4.12) \qquad \bar{\eta}'_\ell(\alpha) \le \bar{\eta}_\ell(\alpha) + \nu(\alpha) + [c_0 - c'_0]^+$$

and

$$(4.13) \qquad \bar{\eta}'_u(\alpha) \le \bar{\eta}_u(\alpha) + [c'_0 - c_0]^+.$$

To see why the claim is true, first note that since $\nu$, $\bar{\eta}_\ell$ and $\bar{\eta}_u$ are nondecreasing, it is clear from (4.10) that if $\bar{\eta}'_\ell$ is continuous at $\alpha$, then (4.12) holds. Likewise, if $\bar{\eta}'_u$ is continuous at $\alpha$, then (4.11) implies that (4.13) is satisfied. Now suppose $\bar{\eta}'_\ell(\alpha) - \bar{\eta}'_\ell(\alpha-) > 0$. Then the complementarity conditions in (1.6) show that $\bar{\phi}'(\alpha) = 0$ and $\bar{\eta}'_u(\alpha-) = \bar{\eta}'_u(\alpha)$. Hence, (4.9) applied to $\bar{\eta}'_\ell$ implies that

$$\bar{\eta}'_\ell(\alpha) = \bar{\eta}'_\ell(\alpha-) - \bar{\phi}'(\alpha-) - \psi'(\alpha) + \psi'(\alpha-).$$

Making the further substitutions $\bar{\eta}'_\ell(\alpha-) - \bar{\phi}'(\alpha-) + \psi'(\alpha-) = -c'_0 + \bar{\eta}'_u(\alpha-)$, $\psi = \psi' + \nu$ and then $\bar{\eta}_u(\alpha-) = c_0 + \psi(\alpha-) + \bar{\eta}_\ell(\alpha-) - \bar{\phi}(\alpha-)$ into the last display, we obtain

$$\begin{aligned}
\bar{\eta}'_\ell(\alpha) &= -c'_0 + \bar{\eta}'_u(\alpha-) - \psi'(\alpha) \\
&= -c'_0 + \bar{\eta}'_u(\alpha-) - \psi(\alpha) + \nu(\alpha) \\
&= -c'_0 + c_0 + \psi(\alpha-) + \bar{\eta}_\ell(\alpha-) - \bar{\phi}(\alpha-) - \psi(\alpha) + \nu(\alpha) \\
&\quad + \bar{\eta}'_u(\alpha-) - \bar{\eta}_u(\alpha-).
\end{aligned}$$

Taking limits as $s \uparrow \alpha$ in (4.11) yields the inequality $\bar{\eta}'_u(\alpha-) - \bar{\eta}_u(\alpha-) \leq [c'_0 - c_0]^+$. When substituted into the last display, this shows that

$$\begin{aligned}
\bar{\eta}'_\ell(\alpha) &\leq -c'_0 + c_0 + \psi(\alpha-) + \bar{\eta}_\ell(\alpha-) \\
&\quad - \bar{\phi}(\alpha-) - \psi(\alpha) + \nu(\alpha) + [c'_0 - c_0]^+ \\
&= \psi(\alpha-) + \bar{\eta}_\ell(\alpha-) - \bar{\phi}(\alpha-) - \psi(\alpha) + \nu(\alpha) + [c_0 - c'_0]^+.
\end{aligned} \tag{4.14}$$

Since $\bar{\eta}_u(\alpha) - \bar{\eta}_u(\alpha-) \geq 0$, (4.9) implies that

$$\begin{aligned}
\bar{\eta}_\ell(\alpha) &= \bar{\eta}_\ell(\alpha-) + [-\bar{\phi}(\alpha-) - \psi(\alpha) + \psi(\alpha-) \\
&\quad + \bar{\eta}_u(\alpha) - \bar{\eta}_u(\alpha-)]^+ \\
&\geq \bar{\eta}_\ell(\alpha-) - \bar{\phi}(\alpha-) - \psi(\alpha) + \psi(\alpha-).
\end{aligned}$$

When substituted into (4.14) this yields (4.12). The proof of the remaining fact that (4.13) continues to hold even if $\bar{\eta}'_u(\alpha) - \bar{\eta}'_u(\alpha-) > 0$ is exactly analogous and is thus omitted.

Having established (4.12) and (4.13), we note from the definition of $\alpha$ that there must exist a sequence $\{s_n\}$ with $s_n \downarrow 0$ as $n \to \infty$ such that one of the following two cases holds:

(4.15) (i) $\quad \bar{\eta}'_\ell(\alpha + s_n) > \bar{\eta}_\ell(\alpha + s_n) + \nu(\alpha + s_n) + [c_0 - c'_0]^+ \qquad \forall n \in \mathbb{N};$

(4.16) (ii) $\quad \bar{\eta}'_u(\alpha + s_n) > \bar{\eta}_u(\alpha + s_n) + [c'_0 - c_0]^+ \qquad \forall n \in \mathbb{N}.$

First, suppose that case (i) holds. Then due to (4.15), the fact that $s_n \downarrow 0$ and the right continuity of $\bar{\eta}'_\ell, \bar{\eta}_\ell$ and $\nu$, it follows that $\bar{\eta}'_\ell(\alpha) \geq \bar{\eta}_\ell(\alpha) + \nu(\alpha) + [c_0 - c'_0]^+$. When combined with (4.12), this yields the equality

$$\bar{\eta}'_\ell(\alpha) = \bar{\eta}_\ell(\alpha) + \nu(\alpha) + [c_0 - c'_0]^+. \tag{4.17}$$

We now show that in this case $\bar{\phi}(\alpha) = \bar{\phi}'(\alpha) = 0$. First, combining (4.17), (4.15) and the fact that $\bar{\eta}_\ell + \nu$ is nondecreasing, we have $\bar{\eta}'_\ell(\alpha + s_n) > \bar{\eta}'_\ell(\alpha)$ for every $n \in \mathbb{N}$. Since $s_n \downarrow 0$, the first complementarity condition in (1.6) ensures that $\bar{\phi}'(\alpha) = 0$. Along with (4.13), (4.17) and the relations $\bar{\phi}'(\alpha) = 0$ and $\psi = \psi' + \nu$, this implies that

$$\begin{aligned}
\bar{\phi}(\alpha) &= \bar{\phi}(\alpha) - \bar{\phi}'(\alpha) \\
&= c_0 - c'_0 + \nu(\alpha) + \bar{\eta}_\ell(\alpha) - \bar{\eta}'_\ell(\alpha) + \bar{\eta}'_u(\alpha) - \bar{\eta}_u(\alpha) \\
&\leq c_0 - c'_0 - [c_0 - c'_0]^+ + [c'_0 - c_0]^+ \\
&= 0.
\end{aligned}$$

Since $\bar{\phi} \in [0, a]$, this implies $\bar{\phi}(\alpha) = 0$.

The right continuity of $\bar{\phi}$ and $\bar{\phi}'$ then ensures the existence of $\varepsilon > 0$ such that for every $s \in [0, \varepsilon]$, $\bar{\phi}(\alpha + s) < a$ and $\bar{\phi}'(\alpha + s) < a$. Hence, due to the complementarity conditions (1.6), property (4.4) and the definitions of $\Gamma_0$ and $\Gamma_{0,a}$, for $s \in [0, \varepsilon]$ we can write

$$\bar{\phi}(\alpha + s) = \Gamma_0(\bar{\phi}(\alpha) + T_\alpha \psi)(s) = \Gamma_0(T_\alpha \psi)(s);$$
$$\bar{\phi}'(\alpha + s) = \Gamma_0(\bar{\phi}'(\alpha) + T_\alpha \psi')(s) = \Gamma_0(T_\alpha \psi')(s);$$
$$T_\alpha \bar{\eta}_\ell(s) = T_\alpha \bar{\eta}(s) = \Gamma_0(T_\alpha \psi)(s) - T_\alpha \psi(s);$$
$$T_\alpha \bar{\eta}_\ell'(s) = T_\alpha \bar{\eta}'(s) = \Gamma_0(T_\alpha \psi')(s) - T_\alpha \psi'(s).$$

Since $T_\alpha \psi = T_\alpha \psi' + T_\alpha \nu$ and $\bar{\phi}(\alpha) = \bar{\phi}'(\alpha) = 0$, property 1 of Lemma 4.1 (replacing $c_0$ and $c_0'$ by 0 and $\psi'$ and $\psi$ by $T_\alpha \psi'$ and $T_\alpha \psi$, resp.) shows that for every $s \in [0, \varepsilon]$,

$$\bar{\eta}_\ell'(\alpha + s) - \bar{\eta}_\ell'(\alpha) = T_\alpha \bar{\eta}_\ell'(s) \leq T_\alpha \bar{\eta}_\ell(s) + T_\alpha \nu(s)$$
$$= \bar{\eta}_\ell(\alpha + s) - \bar{\eta}_\ell(\alpha) + \nu(\alpha + s) - \nu(\alpha).$$

When combined with (4.17) this yields the inequality

$$\bar{\eta}_\ell'(\alpha + s) \leq \bar{\eta}_\ell(\alpha + s) + \nu(\alpha + s) + [c_0 - c_0']^+ \qquad \text{for } s \in [0, \varepsilon],$$

which contradicts (4.15) and so case (i) does not hold.

Thus we have shown that there does not exist any sequence $\{s_n\}$ with $s_n \downarrow 0$ that satisfies (4.15). Together with (4.10) and (4.12), this means that there must exist $\delta > 0$ such that

$$\bar{\eta}_\ell'(s) \leq \bar{\eta}_\ell(s) + \nu(s) + [c_0 - c_0']^+ \qquad \text{for } s \in [0, \alpha + \delta].$$

Combining (1.9) with the above inequality we then obtain for $t \in [0, \alpha + \delta]$,

$$\bar{\eta}_u'(t) = \sup_{s \in [0,t]} [c_0' + \psi'(s) + \bar{\eta}_\ell'(s) - a]^+$$
$$= \sup_{s \in [0,t]} [c_0' + \psi(s) - \nu(s) + \bar{\eta}_\ell'(s) - a]^+$$
$$\leq \sup_{s \in [0,t]} [c_0 + \psi(s) + \bar{\eta}_\ell(s) - a + c_0' - c_0 + [c_0 - c_0']^+]^+$$
$$\leq \sup_{s \in [0,t]} [c_0 + \psi(s) + \bar{\eta}_\ell(s) - a]^+ + [c_0' - c_0]^+.$$
$$= \bar{\eta}_u(t) + [c_0' - c_0]^+.$$

However this contradicts (4.16) and so we conclude that neither case (i) nor case (ii) holds, which in turn contradicts the fact that $\alpha < \infty$. Thus $\alpha = \infty$ or, in other words, the second inequality in property 1 and the first equality in property 2 of the theorem are satisfied.

Applying the result just proved above with $\psi, \psi', c_0, c_0'$ replaced by $-\psi', -\psi$, $a - c_0', a - c_0$ respectively, and invoking (1.10), it follows that $\beta = \infty$, where

$$\beta \doteq \inf\{t > 0 : \bar{\eta}_u'(t) + \nu(t) + [c_0 - c_0']^+ < \bar{\eta}_u(t) \text{ or } \bar{\eta}_\ell'(t) + [c_0' - c_0]^+ < \bar{\eta}_\ell(t)\}.$$

This completes the proof of the first two properties of the theorem. The third and fourth properties are the content of Remark 4.3. □

## APPENDIX: TRANSFORMATION OF GRAPH PARAMETRIZATIONS UNDER $\Lambda_A$

Given $\phi \in \mathcal{D}[0,T]$, recall the definition of the graph $G_\phi$ given in (3.41) and the set $\Pi(\phi)$ of parametric representations of $G_\phi$, as defined immediately after (3.41). The following result is used in the proof of Corollary 1.6.

LEMMA A.1. *Let $\phi \in \mathcal{D}[0,T]$ be given. For $(r,g) \in \Pi(\phi)$, we have $(r, \Lambda_a(g)) \in \Pi(\Lambda_a(\phi))$.*

PROOF. Since the mapping $(r,g)$ is continuous, by Proposition 1.3 the map $(r, \Lambda_a(g))$ is also continuous. We will show that for every $s \in [0,1]$, $(r(s), \Lambda_a(g)(s)) \in G_{\Lambda_a(\phi)}$. Fix $t \in [0,T]$. We consider two cases.

*Case* 1. $\phi(t) = \phi(t-)$.
Consider $s \in [0,1]$ such that $r(s) = t$. We want to show that

$$\Lambda_a(g)(s) = \Lambda_a(\phi)(t), \tag{A.1}$$

which clearly implies $(r(s), \Lambda_a(g)(s)) \in G_{\Lambda_a(\phi)}$. In the case under consideration,

$$g(s) = \phi(t) \tag{A.2}$$

and (A.1) is equivalent to

$$\sup_{s' \in [0,s]} \left[ (g(s') - a)^+ \wedge \inf_{s'' \in [s',s]} g(s'') \right] \tag{A.3}$$
$$= \sup_{t' \in [0,t]} \left[ (\phi(t') - a)^+ \wedge \inf_{t'' \in [t',t]} \phi(t'') \right].$$

The inequality

$$\sup_{s' \in [0,s]} \left[ (g(s') - a)^+ \wedge \inf_{s'' \in [s',s]} g(s'') \right] \tag{A.4}$$
$$\geq \sup_{t' \in [0,t]} \left[ (\phi(t') - a)^+ \wedge \inf_{t'' \in [t',t]} \phi(t'') \right]$$

follows from (A.2) and the monotonicity of $(r,g)$, together with the fact that the graph of $(r(s'), g(s'))$, $s' \in [0,s]$, consists of the graph of $(t', \phi(t'))$, $t' \in [0,t]$, and the vertical segments $\{t'\} \times [\phi(t'-) \wedge \phi(t'), \phi(t'-) \vee \phi(t')]$, $t' \in [0,t]$. To prove the opposite inequality, let $s_0 \in [0,s]$ attain the supremum on the left-hand side of (A.3). Let $t_0 = r(s_0)$ and let $[b, c] = r^{-1}(t_0)$. We want to show that $s_0$ may be chosen to be either $b$ or $c$ (in other words, that the supremum is attained at one of the endpoints of $[b, c]$). This is obvious if $\phi(t_0) = \phi(t_0-)$, since then $g \equiv \phi(t_0)$ on $[b, c]$. If $\phi(t_0-) < \phi(t_0)$, then by the case assumption, $t_0 < t$ and $s_0 \leq c < s$. In this case, $g$ increases on $[b, c]$ and the supremum on the left-hand side of (A.3) is attained at $c$. Thus, if $\phi(t_0-) \leq \phi(t_0)$, we have

$$\sup_{s' \in [0,s]} \left[ (g(s') - a)^+ \wedge \inf_{s'' \in [s',s]} g(s'') \right]$$

$$= (g(c) - a)^+ \wedge \inf_{s'' \in [c,s]} g(s'')$$

$$= (\phi(t_0) - a)^+ \wedge \inf_{t'' \in [t_0,t]} \phi(t'')$$

$$\leq \sup_{t' \in [0,t]} \left[ (\phi(t') - a)^+ \wedge \inf_{t'' \in [t',t]} \phi(t'') \right].$$

On the other hand, if $\phi(t_0-) > \phi(t_0)$, we again have $t_0 < t$ and $s_0 \leq c < s$, but now $g$ decreases on $[b, c]$ and the supremum on the left-hand side of (A.3) is attained at $b$. In this case

$$\sup_{s' \in [0,s]} \left[ (g(s') - a)^+ \wedge \inf_{s'' \in [s',s]} g(s'') \right]$$

$$= (g(b) - a)^+ \wedge \inf_{s'' \in [b,s]} g(s'')$$

$$= (\phi(t_0-) - a)^+ \wedge \phi(t_0-) \wedge \inf_{t'' \in [t_0,t]} \phi(t'')$$

$$\leq \sup_{t' \in [0,t]} \left[ (\phi(t') - a)^+ \wedge \inf_{t'' \in [t',t]} \phi(t'') \right].$$

Thus, regardless of the relationship between $\phi(t_0)$ and $\phi(t_0-)$, (A.3) holds.

*Case* 2. $\phi(t) \neq \phi(t-)$.

Let $[b, c] = r^{-1}(t)$, $\phi' = \phi - (\phi(t) - \phi(t-))\mathbb{I}_{[t,T]}$, $g'(s) = g(s) - (g(s \wedge c) - g(s \wedge b))$. Then $g' = g$ on $[0, b]$, $\phi' = \phi$ on $[0, t)$ and $\phi'(t) = \phi(t-)$. This in turn shows that $\Lambda_a(g')(b) = \Lambda_a(g)(b)$, $\Lambda_a(\phi')(t) = \Lambda_a(\phi)(t-)$ and $(r, g') \in \Pi(\phi')$ on $[0, t]$. Since $\phi'(t) = \phi'(t-)$, we can apply (A.1) to conclude that $\Lambda_a(g)(b) = \Lambda_a(g')(b) = \Lambda_a(\phi')(t) = \Lambda_a(\phi)(t-)$. For $t' > t$ such that $\phi(t') = \phi(t'-)$ and $s' \in [0,1]$ such that $r(s') = t'$ we have, again by (A.1), $\Lambda_a(g)(s') = \Lambda_a(\phi)(t')$. Taking $t' \downarrow t$, we get $\Lambda_a(g)(c) = \Lambda_a(\phi)(t)$. Finally, $\Lambda_a(g)(s)$ moves continuously and monotonically from $\Lambda_a(g)(b)$ to $\Lambda_a(g)(c)$ as $s$ increases over $[b, c]$. Hence, for $s \in [b, c]$, $(r(s), \Lambda_a(g)(s)) = (t, \Lambda_a(g)(s)) \in G_{\Lambda_a(\phi)}$.

*Conclusion.* We have shown that the map $(r, \Lambda_a(g))$ takes values in $G_{\Lambda_a(\phi)}$. If $\Lambda_a(\phi)$ is discontinuous at $t \in (0, T]$, then $\phi$ is also discontinuous at $t$. The Case 2 analysis shows that when $\phi$ is discontinuous at $t$, the function $\Lambda_a(g)$ traverses the vertical segment $\{t\} \times [\phi(t-) \wedge \phi(t), \phi(t-) \vee \phi(t)]$ in the direction from $\phi(t-)$ to $\phi(t)$, which means that $(r, g)$ is nondecreasing in the order relation on the graph of $G_{\Lambda_a(\phi)}$ on the interval $r^{-1}(t)$. For values of $t$ for which $\Gamma_a(\phi)$ is continuous, we use the fact $r$ is nondecreasing to again conclude that $(r, g)$ is nondecreasing. $\square$

**Acknowledgment.** The authors would like to thank Søren Asmussen for bringing the paper [5] to their attention.

L. KRUK  
DEPARTMENT OF MATHEMATICS  
MARIA CURIE-SKLODOWSKA UNIVERSITY  
LUBLIN  
POLAND  
E-MAIL: lkruk@hektor.umcs.lublin.pl

J. LEHOCZKY  
DEPARTMENT OF STATISTICS  
CARNEGIE MELLON UNIVERSITY  
PITTSBURGH, PENNSYLVANIA 15213  
USA  
E-MAIL: jpl@stat.cmu.edu

K. RAMANAN  
S. SHREVE  
DEPARTMENT OF MATHEMATICAL SCIENCES  
CARNEGIE MELLON UNIVERSITY  
PITTSBURGH, PENNSYLVANIA 15213  
USA  
E-MAIL: kramanan@math.cmu.edu  
       shreve@andrew.cmu.edu